\documentclass[12pt]{amsart}

\usepackage{amssymb,amscd,amsthm, verbatim,amsmath,color,fancyhdr, mathrsfs}
\usepackage{graphicx}
\usepackage{turnstile}

\usepackage[letterpaper, left=2.5cm, right=2.5cm, top=2.5cm,
bottom=2.5cm,dvips]{geometry}

\newtheorem{thm}{Theorem}

\newtheorem*{lem*}{Lemma}

\title{UNIQUENESS OF A THREE-DIMENSIONAL STOCHASTIC DIFFERENTIAL EQUATION}

\author[Carl Mueller]{Carl Mueller \textsuperscript{1}} 
\thanks{$^1$Supported by a Simons grant.}
\address{Carl Mueller: Dept. of Mathematics
\\University of Rochester
\\Rochester, NY  14627}
\email{carl.e.mueller@rochester.edu}
\author[Giang Truong]{Giang Truong}
\address{Giang Truong: Dept. of Mathematics
\\University of Rochester
\\Rochester, NY  14627}
\email{gtruong@u.rochester.edu}

\date{\today}

\begin{document}
\maketitle
\begin{abstract}
In order to extend the study of uniqueness property of multi-dimensional systems of stochastic differential equations, in this paper, we look at the following three-dimensional system of equations, of which the two-dimensional case was well-studied before: $dX_t=Y_tdt\quad, dY_t=Z_tdt,\quad dZ_t=|X_t|^{\alpha}dB_t$. We proved that if $(X_0,Y_0,Z_0)\neq(0,0,0)$, and $\frac{3}{4}<\alpha<1$, then the system of equations has a unique solution in the strong sense.
\end{abstract}.
\section{Introduction and main results}
 The uniqueness of ordinary differential equations (ODE) has been extensively studied, see for example \cite{ODE}. In
 particular, if $F(u)$ is Lipschitz continuous, then $$u'(t)=F\left(u(t)\right),\quad u(0)=u_0$$ has a unique
solution for all $t\geq 0$. In the case above, $F$, $u(t)$, and $u_0$ take value in $\mathbf{R}^d$, $d\geq1$. The stochastic differential equation (SDE) realm, on the contrary,
has different criteria for uniqueness of solutions, see for example \cite{SDE}. One of the most well-known results
regarding strong uniqueness of SDE is due to Yamada and
Watanabe's paper \cite{ST}. The result states that if
$f(x)$ is locally H\"older continuous with index $\alpha\in[\frac{1}{2},1]$ and with linear growth, then $$dX=f(X)dW,\quad   X_0=x_0$$ has 
a unique strong solution for all time $t\geq 0.$ Yamada and Watanabe's theory essentially focuses on
one-dimensional SDE.

One motivation for studying higher-dimensional SDE comes from the wave equation:
\begin{equation}\nonumber
\begin{array}{rcl}
    \partial_t^2u & =&\Delta u \\
  u(0,x)&=&u_0(x)\\
  \partial_{t}u(0,x)&=&u_1(x).
\end{array}
\end{equation}
 In this equation, we have $$ \partial_t^2u  =\partial_x^2u=\Delta u.$$ If we let $$ v=\partial_{t} u,$$
 then we can rewrite the wave equation as the following system of equations:
\begin{equation}\nonumber
\begin{array}{rcl}
  \partial_{t} u &= & v \\
    \partial_{t} v&=&\Delta u.
\end{array}
\end{equation}

The original wave equation includes no noise. However, many physical systems are affected by noise. Hence, a modification
of the wave equation which includes white noise is also studied:
\begin{equation}
\label{wave_eq}
    \begin{array}{rcl}
         \partial_{t}^2u & =&\Delta u +f(u)\dot{W}\\
  u(0,x)&=&u_0(x)\\
  \partial_{t}u(0,x)&=&u_1(x).
\end{array}
\end{equation}
    Note: $x\in\mathbf{R} $ and $\dot{W}=\dot{W}(t,x)$ is white noise.
    
One well-known point is that Lipschitz continuity is sufficient for the uniqueness of SDE. Thus, many
mathematicians have studied whether H\"older continuity can still ensure the uniqueness property of SDE.
In Gomez, Lee, Mueller, Neuman, and Salins\cite{GLMNS}, the uniqueness property of the following two-dimensional model of SDE was studied:
\begin{equation}
\label{2-dimensions}
\begin{array}{rcl}
dX&=&Ydt\\
dY&=&|X|^{\alpha}dB\\
(X_0,Y_0)&=&(x_0,y_0).
\end{array}
\end{equation}
The results focused on $f(x)=|x|^{\alpha}$ since it is a prototype of an equation with H\"older continuous coeffients.
Moreover, \eqref{2-dimensions} is a version of \eqref{wave_eq} when we drop the dependence on $x$, which allows us to study the modified wave equation
with more simplicity. We can notice that if we take the differential $dY$ of the first derivative of $X$, which is $Y$
in the system of equations, it resembles
the second derivative in time in the stochastic wave equation. 

It has been proven in \cite{GLMNS} that if $\alpha>\frac{1}{2}$ and $(x_0,y_0)\neq(0,0)$, then \eqref{2-dimensions} has a
unique solution in the strong sense up to the time $\tau$ at which $(X_t,Y_t)$ first hits the
origin $(0,0).$

Since much is still unknown about higher-dimensional SDE, we wish to continue the study of uniqueness property
of \eqref{2-dimensions} in the three-dimensional case, which is: 

\begin{equation}
\label{3-dimensions}
    \begin{array}{rcl}
    dX&=&Ydt \\
   dY&=&Zdt\\
   dZ&=&|X|^ {\alpha}dB\\ 
   (X_0,Y_0,Z_0)&=&(x_0,y_0,z_0).
\end{array}
\end{equation}
  
\begin{thm}
If \(\frac{3}{4}<\alpha<1\) and $(X_0,Y_0,Z_0)\neq0$, then \eqref{3-dimensions} has a unique strong solution, up to the time \(\tau\) at
which the solution $(X_t,Y_t,Z_t)$ first hits the value (0,0,0) or blows up.
\end{thm}

Moreover, we say the solution $(X_t,Y_t,Z_t)$ blows up in finite
time, with positive probability, if there is a random time
$\tau<\infty$ such that $$P\left(\lim_{t\uparrow\tau}|(X_t,Y_t,Z_t)|_{l^\infty}=\infty\right)>0.$$
\section{proof of theorem 1}
Let $(X_t^i,Y_t^i,Z_t^i),i=1,2$ be two solutions to \eqref{3-dimensions} with $(x_0,y_0,z_0)\neq(0,0,0),$ in other words, $(X_t^i,Y_t^i,Z_t^i), i=1,2$ have the same initial condition $(X_0,Y_0,Z_0)\neq (0,0,0)$. Note that since $(X_t^i,Y_t^i,Z_t^i),i=1,2$ have the same initial conditions, from now, we use $X_0,$ $Y_0,$ $Z_0$ instead of $X_0^{i,n},$ $Y_0^{i,n},$ $Z_0^{i,n}.$

Let $\tau$ be the first time $t$ that either $(X_t^1,Y_t^1,Z_t^1)$ or $(X_t^2,Y_t^2,Z_t^2)$
hits the origin $(0,0,0)$ or blows up. We let $\tau$ be infinity if there is no such time.

First of all, if  $X_0\neq0,$ then \eqref{3-dimensions} will have Lipschitz continuity up to the time that $X_t=0$, thus enable uniqueness
to hold. Suppose after a certain amount of time, $X_t$ hits zero, where Lipschitz continuity no longer holds, then due to strong Markov property, we begin the process again with $X_0=0.$ So our goal is to
prove pathwise uniqueness between excursion of X up to the time $\tau$ starting from $X_0=0.$

For any fixed $n$, let $\tau_n$ be the first time that either
$$|(X_t^1,Y_t^1,Z_t^1)|_{l^{\infty}}\wedge |(X_t^2,Y_t^2,Z_t^2)|_{l^{\infty}}\leq 2^{-n}$$
or $$|(X_t^1,Y_t^1,Z_t^1)|_{l^{\infty}}\vee |(X_t^2,Y_t^2,Z_t^2)|_{l^{\infty}}\geq 2^{n}.$$
\vspace{0.3cm}
In which: 

The infinity norm (also known as the $L_\infty$-norm, $l_\infty$-norm, max norm, or uniform norm) of a vector $\vec{v}$ is denoted $|\vec{v}|_{l_\infty}$ and is defined as the maximum of the absolute values of its components:
$$|\vec{v}|_{l_\infty}=\max{ |v_i| : i=1,2,..n},$$

and 
\begin{equation}\nonumber
\begin{array}{rcl}
   a\wedge b  & =& \min(a,b)\\
   a\vee b  &=& \max(a,b).
\end{array}
\end{equation}
If there is no such time, we let $\tau_n$ be infinity. Note that: $\lim_{n \uparrow\infty}(\tau_{n})=\tau$.

Now, for each fixed $n$, we will show uniqueness up to time $\tau_n$ in the following system of equations:
\begin{equation}
\label{3-dim-i,n}
 \begin{array}{rcl}
 dX_t^{i,n}&=&Y_t^{i,n}dt \\
   dY_t^{i,n}&=&Z_t^{i,n}dt\\
   dZ_t^{i,n}&=&|X_t^{i,n}|^ \alpha\mathbf{1}_{[0,\tau_n]}(t)dB_t\\ 
   (X_0,Y_0,Z_0)&=&(x_0,y_0,z_0),
   \end{array}
\end{equation}
 In other words, after the time $\tau_n,$ $dZ_t^{i,n}=0$, which makes $Z_t^{i,n}$ become constant. Specifically, given
 $m,n \in N,$ we need: $$(X_t^n, Y_t^n, Z_t^n)= (X_t^m, Y_t^m, Z_t^m)$$ for all $t\leq\tau_n\wedge\tau_m.$
 
Now, before continuing the proof of uniqueness, by contradiction method, we show that the times that $X_t$ hits zero do not accumulate before the time $\tau_n$, almost surely.

For each $n$, let $A_n$ be the event on which the times that $X_t^{i,n}=0$, $i=1$ or $i=2$, accumulate before $\tau_n$, and assume
$P(A_n)>0.$ Then, on $A_n$, suppose $\sigma_n$ is an accumulation point of the
times $t$ at which $X_t^{i,n}=0,$ i.e, there
exits a sequence of time $\rho_{1,n}<\rho_{2,n}<\ldots$ that converges to
$\sigma_n,$ and $X_{\rho_{k,n}}^{i,n}=0.$ Hence, on $A_n,$ $\lim_{k\to\infty}\rho_{k,n}=\sigma_n.$

We have $X_t^{i,n}$ is almost surely continuous, and that $X_{\rho_{k,n}}^{i,n}=0$ on $A_n$, so: $$\lim_{\rho_{k,n}\to\sigma_n}X_{\rho_{k,n}}^{i,n}=X_{\sigma_n}^{i,n}=0$$ on $A_n.$

Note that $dX_t^{i,n}=Y_t^{i,n}dt,$ and $Y_t^{i,n}$ is almost surely continuous. So if $Y_{\sigma_n}^{i,n}\neq0$ on $A_n,$ then there exists a
random interval $[\sigma_n(\omega)-\epsilon(\omega),\sigma_n(\omega)]$ of positive length for which $X_t^{i,n}\neq0$ on $[\sigma_n(\omega)-\epsilon(\omega),\sigma_n(\omega)].$ This contradicts the hypothesis of $\rho_{k,n}$ converging to $\sigma_n.$

If $Y_{\sigma_n}^{i,n}=0$, there are two cases, which are $$\sigma_n\geq\tau_n,$$ or $$\sigma_n<\tau_n.$$
If $\sigma_n\geq\tau_n,$ then it means that the times at which $X_t^{i,n}$ hit zero do not accumulate before $\tau_n$.

If $\sigma_n<\tau_n,$ then with $X_{\sigma_n}^{i,n}=0$ and $Y_{\sigma_n}^{i,n}=0$, we have $|Z_{\sigma_n}^{i,n}|>2^{-n}$.
Now suppose $Z_{\sigma_n^{i,n}}>2^{-n},$ as the case $Z_{\sigma}^{i,n}<2^{-n}$ is approached similarly due to symmetry.

So, if $Z_{\sigma_n}^{i,n}>2^{-n},$ since $Z_t$ is almost surely continuous, there exists a time interval $[\sigma_{n}(\omega)-\epsilon'(\omega),\sigma_{n}(\omega)]$ on which
$Z_t^{i,n}>\frac{2^{-n}}{2}.$ Hence, for all $t\in[\sigma_n(\omega)-\epsilon'(\omega),\sigma_n(\omega)],$ we almost surely have:

\begin{equation}\nonumber
 \begin{array}{rcl}
   Y_t^{i,n}&=&Y_{\sigma_n}^{i,n} - \int_t^{\sigma_n}Z_s^{i,n}ds  \\[10pt]
     &< & \frac{2^{-n}}{2}(\sigma_n-t).\\[10pt]
 \end{array}   
\end{equation}

Hence, integrating over $X_t^{i,n}$ for $t\in[\sigma_n-\epsilon',\sigma_n],$ we have:
\begin{equation}\nonumber
    \begin{array}{rcl}
       X_t^{i,n}& = & X_{\sigma_n}^{i,n}-\int_t^{\sigma_n}Y_s^{i,n}ds \\[10pt]
        &> & -\int_t^{\sigma_n}Y_s^{i,n}ds\\[10pt]
        &>&\int_t^{\sigma_n}\frac{2^{-n}}{2}(\sigma_n-s)ds\\[10pt]
        &=&\frac{2^{-n}}{2}(\sigma_ns-\frac{s^2}{2})\left. \right\vert_t^{\sigma_n}\\[10pt]
        &=&\frac{2^{-n}}{4}(\sigma_n-t)^2\\[10pt]
        &>&0,
        \end{array}
\end{equation}
as $t<\sigma_n.$ Hence, this contradicts the hypothesis of $\rho_{k,n}$ converging to $\sigma_n$ on $A_n.$
So in conclusion, $P(A_n)=0,$ which means the times at which $X_t$ hits zero does not accumulate before the time $\tau_n.$

One more point we need to address before continuing with the proof of uniqueness is the existence of solutions. This problem is resolved in Theorem 21.7, Theorem 21.8, and Lemma 21.17 of \cite{kallbook},
which prove that for all time $t\geq0$, solutions of multi-dimensional SDE exist with probability one provided the coefficients are continuous and bounded. 

Specifically, in our problem, since $\alpha\in(0,1),$ we have the coefficients of the system of equation \eqref{3-dimensions} are continuous and bounded by $(2^n)^{\alpha}\vee2^n=2^n$ up 
to the $\tau_n$ for each $n$. which satisfies the condition stated in the existence theory in \cite{kallbook}. Hence, existence of solution
holds for all $t\leq\tau_n$ for all n. Therefore, up to the
time $\tau=\sup\tau_n$,  existence of solutions is ensured.

From now, $X_t^{i,n}$ means $X_t^{1,n}$ and $X_t^{2,n}$. Similar interpretation for $Y_t^{i,n}$ and $Z_t^{i,n}$. Back to the proof of uniqueness, we have, with all $t\in[0,\tau_n]:$ $$|Z_t^{i,n}|\vee|Y_t^{i,n}|\leq2^n.$$
So:
\begin{equation}\nonumber
    \begin{array}{rcl}
     |Y_t^{i,n}|& =& |Y_0 +\int_0^t Z_s^{i,n}ds| \\[10pt]
     &\leq&|Y_0|+\int_0^t|Z_s^{i,n}|ds\\[10pt]
      & \leq& 2^n+2^nt.
    \end{array}
\end{equation}
Therefore:
\begin{equation}\nonumber
    \begin{array}{rcl}
     |X_t^{i,n}|&\leq& |X_0+\int_0^tY_s^{i,n}ds|\\[10pt]
         & \leq&\int_0^t|Y_s^{i,n}|ds\\[10pt]
         &\leq&\int_0^t(2^n+2^ns)ds\\[10pt]
         &=&2^n(t+\frac{t^2}{2}).
    \end{array}
\end{equation}
Now let  \begin{equation}
\label{eq:t0}
       t_{0,n}=\frac{2^{-2n}}{2}.
   \end{equation}
Then:
\begin{equation}\nonumber
  \begin{array}{rcl}
      t_{0,n}+\frac{t_{0,n}^2}{2}&= &\frac{2^{-2n}}{2} +\frac{2^{-2n}}{8} \\[10pt]
       &\leq& \frac{2^{-2n}}{2}+\frac{2^{-2n}}{2}\\[10pt]
       &=&2^{-2n}.
  \end{array}  
\end{equation}
So 
\begin{equation}\nonumber
  \begin{array}{rcl}
      2^n(t_{0,n}+\frac{t_{0,n}^2}{2})&\leq &2^n.2^{-2n}  \\[10pt]
      &= & 2^{-n}.
  \end{array}
\end{equation}
Since the quadratic function $2^n(t+\frac{t^2}{2})$ is increasing when $t
\geq0,$ we have:
\begin{equation}\nonumber
   \begin{array}{rcl}
     |X_t^{i,n}|&\leq&2^n(t+\frac{t^2}{2})    \\[10pt]
        &\leq& 2^{-n}
   \end{array} 
\end{equation}
for all $t\in[0,t_{0,n}].$
 
Since $|X_t^{1,n}|$ and $|X _t^{2,n}|$ belong in $ [0,2^{-n}]$ for $t\in [0,t_{0,n}]$, based on the definition of $\tau_n$ above, either
\begin{equation}
\label{case_6}
 |Y_0|\geq2^{-n} 
\end{equation}
or 
\begin {equation}
\label{case_7}
|Z_0|\geq 2^{-n}
\end{equation}
for each fixed $n$. This is due to the fact
that the solutions have the same initial condition $(X_0,Y_0,Z_0)$ and for each time $t\in[0,t_{n,0}]$, either :
\begin{equation}\nonumber
    |Y_t^{i,n}|\geq2^{-n}
\end{equation}
or
\begin{equation}\nonumber
   |Z_t^{i,n}|\geq2^{-n}.
\end{equation}
First, we deal with $Y_0 > 0$ since due to symmetry, we could deal with the case $Y_0<0$ with similar methods.

Now, with $Y_0 > 0$, we look at other subcases based on $Z_0^{i,n}.$

\textbf{Case I:} $Y_0>0,$ $|Z_0|\leq 2^{-n}$ .

If $|Z_0|\leq 2^{-n}$, then \eqref{case_6} takes place. We are looking at the case $Y_0>0,$ thus $Y_0>2^{-n}.$  Also, note that $dZ_t^{i,n}=0\quad\forall t>\tau_n$, and $|Z_t^{i,n}|\leq 2^{n}$ for $\forall t\in[0,\tau_n].$  Hence, $|Z_t^{i,n}|\leq2^n\quad\forall t$, which means $Z_t^{i,n}\geq-2^n \quad \forall t.$
Next, we have:
\begin{equation}\nonumber
    \begin{array}{rcl}
    Y_t^{i,n}&=&Y_0 +\int_0^t Z_sds\\[10pt]
    &\geq&2^{-n} - \int_0^t 2^nds\\[10pt]
    &=&2^{-n}-2^nt.\\
    \end{array}
\end{equation}

If  $$0<t<t_{0,n}=\frac{2^{-2n}}{2},$$ where $t_{0,n}$ is defined as in \eqref{eq:t0}, then:

    $$ 2^nt< \frac{2^{-n}}{2} $$
     thus:
      $$2^{-n}-2^nt> \frac{2^{-n}}{2}$$
$ \forall t\in[0,t_{0,n}].$
In other words, $Y_t^{i,n}>\frac{2^{-n}}{2}$ for $t\in[0,t_{0,n}].$
So, for all $t\in[0,t_{0,n}],$ we have:$$ Y_t^{i,n}\geq\frac{2^{-n}}{2}.$$ Hence:
\begin{equation}
\label{8}
    \begin{array}{rcl}
    X_t^{i,n}& \geq & X_0^{i,n} + \int_{0}^{t}\frac{2^{-n}}{2}ds \\[10pt]
    &\geq& \frac{2^{-n}}{2}t.
    \end{array}
\end{equation}

Furthermore, based on \eqref{8} and $|X_t^{i,n}|\leq 2^{-n}$ for all $t\in[0,t_{0,n}]$, it leads to  $t_{0,n}\leq 2,$ otherwise $ X_t^{i,n}>2^{-n},$ which means that $t>t_{0,n},$ hence a contradiction.

Note that
\begin{equation}\nonumber
\begin{array}{rcl}
 X_t^{i,n}& =&X_0+\int_{0}^{t}Y_s^{i,n}ds \\[10pt]
 Y_s^{i,n}&=    &Y_0+\int_{0}^{s}Z_k^{i,n}dk\\[10pt]
 Z_k^{i,n}&=&Z_0+\int_{0}^{k}|X_r^{i,n}|^{\alpha}\mathbf{1}_{[0,\tau_n]}(t)dB_r.
\end{array}
\end{equation}
Hence:
\begin{equation}\nonumber
    \begin{array}{rcl}
     X_t^{i,n}&  =   &X_0+Y_{0}t+\displaystyle\int_{0}^{t}\int_{0}^{s}\left(Z_0+\int _{0}^{k}|X_r^{i,n}|^{\alpha}\mathbf{1}_{[0,\tau_n]}(t)dB_r\right)dkds   \\[10pt]

         &=& X_0+Y_0t+\displaystyle\int_{0}^{t}\int_{0}^{s} Z_0dkds+\int_{0}^{t}\int_{0}^{s}\int_{0}^{k}|X_r^{i,n}|^{\alpha}\mathbf{1}_{[0,\tau_n]}(t)dB_{r}dkds\\[10pt]
         
         &=&X_0+Y_0t +Z_0\frac{t^2}{2} +\displaystyle\int_{0}^{t}\int_{0}^{s}\int_{0}^{k}|X_r^{i,n}|^{\alpha}\mathbf{1}_{[0,\tau_n]}(t)dB_rdkds.\\
    \end{array}
\end{equation}
 
Hence:
$$(X_t^{1,n}-X_t^{2,n})^2=\left(\int_{0}^{t}\int_{0}^{s}\int_{0}^{k}\left(|X_{r}^{1,n}|^{\alpha}-|X_{r}^{2,n}|^{\alpha}\right)\mathbf{1}_{[0,\tau_n]}(r)dB_rdkds\right)^2.$$
Apply the Cauchy-Schwarz inequality twice:
$$\leq t\int_{0}^{t}
\left(\int_{0}^{s}\int_{0}^{k}\left(|X_{r}^{1,n}|^{\alpha}-|X_{r}^{2,n}|^{\alpha}\right)\mathbf{1}_{[0,\tau_n]}(r)dB_rdk\right)^2ds $$
$$\leq t\int_{0}^{t}s\int_{0}^{s}\left(\int_{0}^{k}\left(|X_{r}^{1,n}|^{\alpha}-|X_{r}^{2,n}|^{\alpha}\right)\mathbf{1}_{[0,\tau_n]}(r)dB_r\right)^2dkds$$
Thus:
$$ E\left[(X_t^{1,n}-X_t^{2,n})^2\right]\leq tE\int_{0}^{s}s\int_{0}^{s}\left(\int_{0}^{k}\left(|X_{r}^{1,n}|^{\alpha}-|X_{r}^{2,n}|^{\alpha}\right)\mathbf{1}_{[0,\tau_n]}(r)dB_r\right)^2dkds.$$
By Ito's isometry:
\begin{equation}  
\begin{array}{rcl}\nonumber
 E\left[(X_t^{1,n}-X_t^{2,n})^2\right] &\leq& t E\displaystyle\int_{0}^{t}s\int_{0}^{s}\int_{0}^{k}\left((|X_{r}^{1,n}|^{\alpha}-|X_{r}^{2,n}|^{\alpha})\mathbf{1}_{[0,\tau_n]}(r)\right)^2drdkds \\[10pt]
 &\leq&  tE \displaystyle\int_{0}^{t}t\int_{0}^{s}\int_{0}^{k}\left(|X_{r}^{1,n}|^{\alpha}-|X_{r}^{2,n}|^{\alpha}\right)^2drdkds \\[10pt]
 &\leq& t^2E \displaystyle\int_{0}^{t}\int_{0}^{s}\int_{0}^{k}\left(|X_{r}^{1,n}|^{\alpha}-|X_{r}^{2,n}|^{\alpha}\right)^2drdkds \\[10pt]
 &\leq& t^2E \displaystyle\int_{0}^{t}\int_{0}^{t}\int_{0}^{t}\left(|X_{r}^{1,n}|^{\alpha}-|X_{r}^{2,n}|^{\alpha}\right)^2drdkds \\[10pt]
 &=&t^4E\displaystyle\int_{0}^{t}\left(|X_{r}^{1,n}|^{\alpha}-|X_{r}^{2,n}|^{\alpha}\right)^2dr.
 
\end{array}
\end{equation}
Now we apply the Mean Value Theorem for the function $f(x)=x^\alpha, 0< \alpha<1$, and $a<b$ :
$$b^\alpha -a^\alpha =\alpha c^{\alpha -1}(b-a)\leq \alpha a^{\alpha -1}(b-a)$$
for $c\in(a,b).$
Then for $r\in[0,t_{0,n}],$ where $t_0$ is determined in \eqref{eq:t0}, we apply \eqref{8}:
\begin{equation}
\label{method}
 \Big||X_{r}^{1,n}|^{\alpha}-|X_{r}^{2,n}|^{\alpha}\Big|\leq \alpha\left(\frac{2^{-n}}{2}r\right)^{\alpha-1}\Big||X_{r}^{1,n}|-|X_{r}^{2,n}|\Big|.
 \end{equation}

Now let $$D_t=E\Big[\left(|X_{t}^{1,n}|-|X_{t}^{2,n}|\right)^2\Big].$$
Since $t_{0,n}\leq2$, we have for all $t\in[0,t_{0,n}],$
$$D_t\leq E\left[(X_{t}^{1,n}-X_{t}^{2,n})^2\right]\leq C_n\int_{0}^{t}r^{2\alpha -2}D_rdr.$$

for some $C_n$ depending on $n$. Since $\alpha > \frac{3}{4},  r^{4\alpha-4}$ is integrable on $ [0,t_{0,n}].$ At this stage we apply Gronwall's lemma:
\begin{lem*}
Let $I$ denote an interval of the real line of the form $[a,\infty)$ or $[a,b]$ or $[a,b)$ with $a<b.$ Let $\beta$ and
$u$ be real-valued continuous functions defined on $I$. If $u$ is differentiable on the interior $I^{o}$ of $I$ (the
interval $I$ without the endpoint $a$ and possibly $b$) and satisfies the differential inequality
$$ u^{'}(t)\leq\beta(t)u(t),\quad t\in I^{o},$$
then $u$ is bounded by the solution of the corresponding differential equation $v^{'}(t)=\beta(t)v(t):$
$$u(t)\leq u(a)\exp\left(\int_{a}^{t}\beta(s)ds\right)$$
for all $t\in I.$
\end{lem*}
\vspace{0.5cm}
 Hence, with $D_0=0$, we have $D_t=0$ for all $t\in[0,t_{0,n}].$ Therefore, \eqref{3-dimensions} has unique strong solution in $[0,t_{0,n}].$
 
Since $\alpha\geq\frac{3}{4}$, we have $2\alpha -2 \geq -1$. Hence, $r^{2\alpha -2}$ is integrable on $[0,t_{0,n}]$
Note: in this case since $X_t\geq 2^{-n},$ $\eta\leq 1.$ Applying Gronwall's lemma, with $D_0=0,$ we have $D_t=0$ for
all $t\in[0,t_{0,n}]$. Therefore, \eqref{3-dimensions} has a unique solution in strong sense up to time $t_{0,n}$. 

Note that for all $t\in[0,t_{0,n}],$ we have $Y_t^{i,n}>\frac{2^{-n}}{2}>0.$ Hence $X_t^{i,n}$ is strictly increasing, which, leads to $X_{t_{0,n}}$ strictly positive. Therefore, by strong Markov property, we have uniqueness till the time $X$ next hits zero.
Now we proceed to the next case.

\textbf{Case II:} $Y_0>0$, $Z_0< -2^{-n}.$ 

Since $Z_0$ starts negative, $Y_t^{i,n}$ decreases for an amount of time. Since $Y_0$ is positive, let say $Y_0=\beta >0.$
Note that $\forall t,$ we have $|Z_t^{i,n}|<2^{n},$ which means $Z_t^{i,n}>-2^n.$ First, we have;
\begin{equation}\nonumber
\begin{array}{rcl}
    Y_t^{i,n}&= & \beta +\int_0^t Z_sds \\[10pt]
    &\geq& \beta-\int_0^t 2^nds\\[10pt]
    &=&\beta-2^nt.\\[10pt]
\end{array}
\end{equation}

Let 
\begin{equation}
\label{t'}
t_{0,n}^{'}=\frac{\beta}{2^{n+1}}.
\end{equation}
If $$0<t<t_{0,n}^{'}=\frac{\beta}{2^{n+1}},$$ then:
$$2^nt< \frac{\beta}{2}$$
thus
   $$\beta-2^nt>\frac{\beta}{2} $$
for $\forall t\in[0,t_{0,n}{'}].$ In other words, $Y_t^{i,n}>\frac{\beta}{2}$ for $\forall t\in[0,t_{0,n}^{'}].$ 

So , $\forall t\in[0,t_{0,n}\wedge t_{0,n}^{'}],$ where $t_{0,n}$ and $t_{0,n}^{'}$ are determined in \eqref{eq:t0} and \eqref{t'} respectively, we have:

  $$ Y_t^{i,n}\geq\frac{\beta}{2}$$
  Hence:
\begin{equation}\nonumber
\begin{array}{rcl}
      X_t^{i,n}&\geq & X_0^{i,n} + \int_{0}^{t}\frac{\beta}{2}ds \\[10pt]
     &\geq& \frac{\beta}{2}t.\\
     \end{array}
 \end{equation}
 
Applying the same method \eqref{method} above, we use Mean Value Theorem  for the new lower bound of $X_t^{i,n}$:
$$\Big||X_{r}^{1,n}|^{\alpha}-|X_{r}^{2,n}|^{\alpha}\Big|\leq \alpha\left(\frac{\beta}{2}r\right)^{\alpha-1}\Big||X_{r}^{1,n}|-|X_{r}^{2,n}|\Big|.$$
Hence:
$$ D_t\leq E\left[(X_{t}^{1,n}-X_{t}^{2,n})^2\right]\leq C_n\int_{0}^{t}r^{2\alpha -2}D_rdr.$$
Again, applying Gronwall's lemma, with $D_0=0$, we have $D_t=0$ for all $t\in[0,t_{0,n}\wedge t_{0,n}^{'}].$ Therefore, \eqref{3-dimensions} has unique strong solution in $[0,t_{0,n}\wedge t_{0,n}^{'}].$
As in the previous cases, we have $Y_t^{i,n}>\frac{\beta}{2}>0$ for all $t\in[t_{0,n}\wedge t_{0,n}^{'}],$ which makes $X_t^{i,n}$ strictly increasing. So $X_{t_{0,n}\wedge t_{0,n}^{'}}$ is strictly positive. Thus, by strong Markov property, we have uniqueness till the next time $X$ hits zero.

\textbf{Case III: }
$Y_0=0$, $Z_0 > 2^{-n}$. 

Now we let $T_n$ be the first time that either $Z_t^{1,n}$ or $Z_t^{2,n}$ hits the value $\frac{2^{-n}}{2}.$ Since both $Z_t^{1,n}$ and $Z_t^{2,n}$ are continuous, we have $T_n>0$ with probability one. So we now prove uniqueness up
to the time $t_{0,n}\wedge T_n,$ where $t_{0,n}$ is defined in \eqref{eq:t0}.

 Then for all $t$ in $[0,t_{0,n}\wedge T_{n}],$ we have: $$Z_t^{i,n}\geq\frac{2^{-n}}{2}, $$ therefore:
\begin{equation}
\label{case_11}
\begin{array}{rcl}
Y_t^{i,n}&\geq& Y_0 + \int_{0}^{t}\frac{2^{-n}}{2}ds\\[10pt]
  &\geq& \frac{2^{-n}}{2}t.
    \end{array}
\end{equation}
(Because $Y_0^{i,n}\geq 0$).

Based on \eqref{case_11} for $t\in[0,t_{0,n}\wedge T_n]:$
\begin{equation}
\label{12}
    \begin{array}{rcl}
      X_t^{i,n}&\geq   & X_0+\int_{0}^{t}\frac{2^{-n}}{2}s ds \\[10pt]
     & \geq&\frac{2^{-n}}{4}t^{2}.
    \end{array}
\end{equation} 
Now we denote:
\begin{equation}\nonumber
   \begin{array}{rcl}
  X_t^{i,n}&=&\tilde{X}_t^{i,n}\\[10pt]
  Y_t^{i,n}&=&\tilde{Y}_t^{i,n}\\[10pt]
  Z_t^{i,n}&=&\tilde{Z}_t^{i,n}\\
   \end{array}
\end{equation}
for $i=1,2$ and for $t\leq\tau_n\wedge\ T_{n}\wedge t_{0,n},$ as $t_{0,n}$ defined as in \eqref{eq:t0} above.

Hence, the following system of equations holds up to the stopping time $\tau_n\wedge T_{n}\wedge t_{0,n}.$
\begin{equation}
\label{tilde-system}
    \begin{array}{rcl}
      d\tilde{X}_t^{i,n}&=   &\tilde{Y}_t^{i,n}dt  \\[10pt]
        d\tilde{Y}_t^{i,n} & =&\tilde{Z}_t^{i,n}dt\\[10pt]
        d\tilde{Z}_t^{i,n}&=&|\tilde{X}_t^{i,n}|^{\alpha}\mathbf{1}_{[0,\tau_n\wedge T_{n}\wedge t_{0,n}]}(t)dB_t,
    \end{array}
\end{equation}
with $(\tilde{X}_0^{i,n},\tilde{Y}_0^{i,n},\tilde{Z}_0^{i,n})=(X_0,Y_0,Z_0),$ for $i=1,2.$ Furthermore, using \eqref{tilde-system}, $\tilde{X}_t^{i,n}, \tilde{Y}_t^{i,n},$ and 
$\tilde{Z}_t^{i,n}$ can be defined for all time.
 
Use Ito's isometry as above with $\tilde{X}_t^{i,n}$, $\tilde{Y}_t^{i,n}$, and $\tilde{Z}_t^{i,n}$:
\begin{equation}  
\begin{array}{rcl}\nonumber
 E\left[(\tilde{X}_t^{1,n}-\tilde{X}_t^{2,n})^2\right] &\leq& t E\displaystyle\int_{0}^{t}s\int_{0}^{s}\int_{0}^{k}\left((|\tilde{X}_r^{1,n}|^{\alpha}-|\tilde{X}_r^{2,n}|^{\alpha})\mathbf{1}_{[0,\tau_n\wedge T_{n}\wedge t_{0,n}]}(r)\right)^2drdkds \\[10pt]
 &\leq&  tE \displaystyle\int_{0}^{t}t\int_{0}^{s}\int_{0}^{k}\left(|\tilde{X}_{r}^{1,n}|^{\alpha}-|\tilde{X}_{r}^{2,n}|^{\alpha}\right)^2\mathbf{1}_{[0,\tau_n\wedge T_{n}\wedge t_{0,n}]}(r)drdkds \\[10pt]
 &\leq& t^2E \displaystyle\int_{0}^{t}\int_{0}^{s}\int_{0}^{k}\left(|\tilde{X}_{r}^{1,n}|^{\alpha}-|\tilde{X}_{r}^{2,n}|^{\alpha}\right)^2\mathbf{1}_{[0,\tau_n\wedge T_{n}\wedge t_{0,n}]}(r)drdkds \\[10pt]
 &\leq& t^2E \displaystyle\int_{0}^{t}\int_{0}^{t}\int_{0}^{t}\left(|\tilde{X}_{r}^{1,n}|^{\alpha}-|\tilde{X}_{r}^{2,n}|^{\alpha}\right)^2\mathbf{1}_{[0,\tau_n\wedge T_{n}\wedge t_{0,n}]}drdkds \\[10pt]
 &=&t^4E\displaystyle\int_{0}^{t}\left(|\tilde{X}_{r}^{1,n}|^{\alpha}-|\tilde{X}_{r}^{2,n}|^{\alpha}\right)^2\mathbf{1}_{[0,\tau_n\wedge T_{n}\wedge t_{0,n}]}dr.
\end{array}
\end{equation}

Using \eqref{12} and Mean Value Theorem, for $r\in[0,\tau_n\wedge T_{n}\wedge t_{0,n}],$ we have:
\begin{equation}\nonumber
\Big||\tilde{X}_{r}^{1,n}|^{\alpha}-|\tilde{X}_{r}^{2,n}|^{\alpha}\Big|\leq \alpha\left(\frac{2^{-n}}{4}r^{2}\right)^{\alpha-1}\Big||\tilde{X}_{r}^{1,n}|-|\tilde{X}_{r}^{2,n}|\Big|,
\end{equation}

Hence:
$$ E\left[(\tilde{X}_t^{1,n}-\tilde{X}_t^{2,n})^2\right]\leq t^4\alpha^2\left(\frac{2^{-n}}{4}\right)^{2(\alpha-1)}E\displaystyle\int_{0}^{t}r^{4(\alpha-1)}(|\tilde{X}_{r}^{1,n}|-|\tilde{X}_{r}^{2,n}|)^2dr.$$

So if  we let $$D_t=E\Big[\left(|\tilde{X}_{t}^{1,n}|-|\tilde{X}_{t}^{2,n}|\right)^2\Big],$$ then: 

$$D_t\leq E\left[(\tilde{X}_{t}^{1,n}-\tilde{X}_{t}^{2,n})^2\right]\leq C_n\int_{0}^{t}r^{4\alpha -4}D_rdr.$$ 
Again, applying Gronwall's lemma, with $D_0=0$, we have $D_t=0.$ 
 Note that at the time $t_{0,n}\wedge T_{n},$ since we have $Z_t^{i,n}>0$ for all $t\in[0,t_{0,n}\wedge T_{n}]$, and also initial $Y_0>0$, it leads to $Y_t^{i,n}>0$ for all
 $t\in[0,t_{0,n}\wedge T_{n}].$ Thus $X_t^{i,n}$ is strictly increasing, which makes $X_{t_{0,n}\wedge
 T_{n}}$ must be strictly greater than zero. Therefore, by strong Markov property, we obtain uniqueness of the process till $X$ next hits zero.
Now we continue with the proof of the theorem with the next case.

\vspace{0.5cm}

For $Y_0 < 0,$ as stated above, is solved similarly due to symmetry.
\vspace{0.5cm}

Lastly, we complete the proof of Theorem 1 with:

\textbf{Case IV:} $Y_0=0.$

If $Y_0=0$, then based on the definition of $\tau_n,$
$|Z_0|>2^{-n}.$ We will first deal with the case $Z_0> 2^{-n},$ and the case $Z_0< 2^{-n}$ is
approached the same way due to symmetry. As in \textbf{ Case III}, let $T_n$ be the first
time that either $Z_t^{1,n}$ or $Z_t^{2,n}$ hits the value $\frac{2^{-n}}{2}.$ Due to the continuity of $Z_t^{1,n}$ and $Z_t^{2,n}$, $T_{n}>0$ with probability one. So with $\forall t\in[0,t_0\wedge T_{n}],$ where $t_0$ is determined in \eqref{eq:t0}, we have:
\begin{equation*}\nonumber
    \begin{array}{rcl}
    Y_t^{i,n}&\geq& Y_0 + \int_{0}^{t}\frac{2^{-n}}{2}ds\\[10pt]
  &=&\frac{2^{-n}}{2}t.\\[10pt]
  \end{array}
\end{equation*}
Then
\begin{equation*}
    \begin{array}{rcl}
        X_t^{i,n}&\geq   & X_0+\int_{0}^{t}\frac{2^{-n}}{2}s ds \\[10pt]
     & \geq&\frac{2^{-n}}{4}t^{2}.

    \end{array}
\end{equation*}
We now apply the same method as in \textbf{Case III} by looking at $\tilde{X}_t^{i,n},\tilde{Y}_t^{i,n},$ and $\tilde{Z}_t^{i,n},$ which are defined as:
\begin{equation}\nonumber
   \begin{array}{rcl}
  X_t^{i,n}&=&\tilde{X}_t^{i,n}\\[10pt]
  Y_t^{i,n}&=&\tilde{Y}_t^{i,n}\\[10pt]
  Z_t^{i,n}&=&\tilde{Z}_t^{i,n}\\
   \end{array}
\end{equation}
for $i=1,2$ and for $t\leq\tau_n\wedge\ T_{n}\wedge t_{0,n},$ as $t_{0,n}$ defined as in \ref{eq:t0} above.

Hence, the following system of equations holds up to the stopping time $\tau_n\wedge T_{n}\wedge t_{0,n}.$
\begin{equation}\nonumber
    \begin{array}{rcl}
      d\tilde{X}_t^{i,n}&=   &\tilde{Y}_t^{i,n}dt  \\[10pt]
        d\tilde{Y}_t^{i,n} & =&\tilde{Z}_t^{i,n}dt\\[10pt]
        d\tilde{Z}_t^{i,n}&=&|\tilde{X}_t^{i,n}|^{\alpha}\mathbf{1}_{[0,\tau_n\wedge T_{n}\wedge t_{0,n}]}(t)dB_t,
    \end{array}
\end{equation}
with $(\tilde{X}_0^{i,n},\tilde{Y}_0^{i,n},\tilde{Z}_0^{i,n})=(X_0,Y_0,Z_0),$ for $i=1,2.$ Furthermore, using \eqref{tilde-system}, $\tilde{X}_t^{i,n}, \tilde{Y}_t^{i,n}, \tilde{Z}_t^{i,n}$ can be defined for all time.

Again, using the same strategy in \textbf{Case III} and Mean Value Theorem, we have:
$$ \Big||\tilde{X}_{r}^{1,n}|^{\alpha}-\tilde{X}_{r}^{2,n}|^{\alpha}\Big|\leq \alpha\left(\frac{2^{-n}}{4}r^{2}\right)^{\alpha-1}\Big||\tilde{X}_{r}^{1,n}|-|\tilde{X}_{r}^{2,n}|\Big|.$$
Hence, if we let $$D_t=E\Big[\left(|\tilde{X}_{t}^{1,n}|-|\tilde{X}_{t}^{2,n}|\right)^2\Big],$$ then :
$$D_t\leq E\left[(\tilde{X}_{t}^{1,n}-\tilde{X}_{t}^{2,n})^2\right]\leq C_n\int_{0}^{t}r
^{4\alpha -4}D_rdr.$$
Hence, using Gronwall's lemma with $D_0=0$ yields $D_t=0$, which completes the proof of Theorem 1.

In this casee, we also have $Y_t^{i,n}>\frac{2^{-n}}{2}t>0$ for all $t\in[0,t_{0,n}\wedge T_{n}]$. Hence $X_t^{i,n}$ is strictly increasing, which yields $X_{t_{0,n}\wedge T_{n}}$ strictly positive. So, by strong Markov property, we have uniqueness up to the time  $X$ next hits zero.

Now with uniqueness proved, we actually can even strengthen the proof by showing that $\tau_n^1=\tau_n^2$ for all $n$, where $\tau_n^1$ and $\tau_n^2$ respectively stand for the stopping time at the critical value for $X_t^1$ and $X_t^2.$ Without loss of generality, suppose $\tau_n^2>\tau_n^1:$
So at the time $\tau_n^1$, $X_t^2$ has not yet reached the critical values, which are $2^{-n}$ or $2^n,$ as stated above. But since we have uniqueness up to $\tau_n^1\wedge\tau_n^2,$ this implies $X_t^1$ has also not reached the critical value at the time $\tau_n^1,$ which is a contradiction to the definition of $\tau_n^1.$ Hence, $\tau_n^1$ and $\tau_n^2$ must be equal.
\vskip.5in

\bibliography{mybibliography}
\bibliographystyle{plain}

%

\end{document}